\def\wR{\mbox{$I\!\!R$}}
\def\wD{\mbox{$I\!\!D$}}

\font\textmsbm= msbm10 scaled 1200
\font\scriptmsbm= msbm7 scaled 1200
\font\scriptscriptmsbm= msbm5 scaled 1200
\newfam\tmwfama
\textfont\tmwfama\textmsbm
\scriptfont\tmwfama\scriptmsbm
\scriptscriptfont\tmwfama\scriptscriptmsbm
\def\bbb{\fam\tmwfama\textmsbm}
\def\bT{{\bbb T}}
\def\span{\mbox{\,\rm span\,}}

\def\Lra{\Leftrightarrow}

\def\Ra{\Rightarrow}

\def\d{\delta}
\def\e{\epsilon}

\def\th{\theta}

\def\s{\sigma}

\def\o{\omega}

\def\es{\emptyset}

\def\sm{\setminus}

\def\sb{\subset}      \def\sbe{\subseteq}
\def\sp{\supset}      \def\spe{\supseteq}

\def\tN{\mathchoice{\mbox{\bf N}}{\mbox{\bf N}}{\mbox{\scriptsize \bf
N}}{\mbox{\tiny\bf N}}}

\def\bydef{\,\lower-.1ex\hbox{\rm :}\!=}

\def\cB{{\cal B}}
\def\cC{{\cal C}}
\def\cD{{\cal D}}
\def\cE{{\cal E}}

\def\cG{{\cal G}}

\def\cK{{\cal K}}
\def\cL{{\cal L}}

\def\cP{{\cal P}}

\def\cS{{\cal S}}

\def\cV{{\cal V}}

\def\BMO{\mbox{\,\rm BMO\,}}

\def\ms{\vskip 0.3cm}
\def\bs{\vskip 0.5cm}
\def\proof{\goodbreak\noindent{\sc Proof. }\nobreak}
\def\prooof{\goodbreak\noindent{\sc Proof of Theorem 1. }\nobreak}
\def\proooof{\goodbreak\noindent{\sc Proof of Theorem 2. }\nobreak}

\def\rem{\vskip3pt\noindent{\sc Remark.} }
\def\rems{\vskip3pt\noindent{\sc Remarks.} }
\def\endrems{\vskip3pt}
\def\endrem{\vskip3pt}

\def\endproof{\par\nobreak\hbox to \hsize{\hfil\vrule width 5pt height
5pt}\goodbreak\vskip 3pt}
\def\endprooof{\par\nobreak\hbox to \hsize{\hfil\vrule width 5pt height
5pt}\goodbreak\vskip 3pt}
\def\endproooof{\par\nobreak\hbox to \hsize{\hfil\vrule width 5pt height
5pt}\goodbreak\vskip 3pt}

\def\se{\par\nobreak\hbox to \hsize{\hfil\vrule width 5pt height
5pt}\goodbreak\vskip 3pt}
\newtheorem{theor}{Theorem}
\newtheorem{prop}{Proposition}

\newtheorem{lemma}{Lemma}

\def\1{\mbox{ \bf 1 }}
\documentstyle[12pt,twoside,titlepage]{article}
\title{Rearrangements of the Haar system which preserve $\BMO$}
\author{by\\
Paul F.X. M\"uller$^*$}
\date{}

\footnotetext
{$^*$ Supported by the Austrian Academy of Sciences, APART Programm}
\footnotetext
{1991 Mathematics Subject}
\footnotetext
{Classification: 42C20, 43A17, 47B38}

\titlepage
\begin{document}
\parindent0.5cm
\maketitle
\newpage
\pagestyle{plain}
\pagenumbering{arabic}

\section{Introduction}

Let $\cD$ denote the collection of dyadic intervals in $[0,1]$. Let $\tau$ be a
rearrangement of the dyadic intervals; i.e., $\tau:\cD\to \cD$ is an injective
map. In this paper we study the induced operator given by the relation
$$Th_I=h_{\tau(I)},$$
where $h_I$ is the $L^\infty$-normalized Haar function.

A function $x$ on $[0,1]$ with Haar expansion $\sum x_Ih_I$, is in dyadic $\BMO$
if
$$||x||=\sup_{J\in\cD}\left(\frac 1{|J|}\sum_{I\sbe
J}x^2_I|I|\right)^{1/2}<\infty.$$
$||T||$ is the operator norm of $T$ on $\BMO$.
We give the following
 geometric characterization of those rearrangements $\tau$ for which
$||T||\,||T^{-1}||<\infty$ respectively $||T||<\infty$: $\tau$ and $\tau^{-1}$ satisfy Property
$\cP$ whenever $T$ is an isomorphism on $\BMO$; $\tau$ satisfies the weak-Property
$\cP$ whenever $T$ is a bounded operator, (see Section 2 for definitions).
We will also see that $T$ is bounded whenever $\tau$ preserves the Carleson
packing condition.
(We say that $\cL\sb \cD$ satisfies the $M$ Carleson condition
if,
$$\sup_{J\in\cD}\frac 1{|J|}\sum_{I\in\cL,I\sbe J}|I|\le M.\leqno(1.1)$$
${\lbrack\!\lbrack}\cL{\rbrack\!\rbrack}$ denotes the infimum of the constants
$M$ that satisfy (1.1).)

The following two theorems summarize the results of this paper:

\begin{theor} For a rearrangement $\tau$ the following are equivalent.

(i) $T:\BMO\to\BMO$ is an isomorphism.

(ii) There exists $M\ge 0$ so that for any $\cL\sb\cD$,\newline $M^{-1}
{\lbrack\!\lbrack}\cL{\rbrack\!\rbrack}\le{\lbrack\!\lbrack}\tau(\cL)
{\rbrack\!\rbrack}\le
M\lbrack\!\lbrack\cL\rbrack\!\rbrack$.

(iii) $\tau$
and $\tau^{-1}$
satisfy Property $\cP$.
\end{theor}

\begin{theor} For a rearrangement $\tau$ the following are equivalent.

(i) $T:\BMO\to \BMO$ is a bounded operator.

(ii) There exists $M\ge 1$ so that for $\cL\sb\cD$, ${\lbrack\!\lbrack}\tau(\cL)
\rbrack\!\rbrack\le M{\lbrack\!\lbrack}\cL{\rbrack\!\rbrack}$.

(iii) $\tau$
satisfies the weak-Property $\cP$.
\end{theor}


We should emphasise the fact that the theorems above concern
rearrangements of the $L^\infty$ normalized Haar system. For, had we
considered rearrangements of the $L^2$ normalized Haar system in $\BMO$, then
the difficulties of the above problems would have disappeared almost entierly.

The first significant results about the behaviour of the Haar system under
rearrangements were due to E.M. Semyonov, see [S1]-[S3]. Under the a-priori
assumption that $|\tau(I)|=|I|$ for $I\in\cD$ E.M. Semyonov found a necessary and
sufficient condition for the boundednes of the induced permutation operator
acting from $L^p$ into $L^q$. The present work is partly motivated by the desire
to eliminate any a-priori assumptions from Semyonov's theorems.

A further motivation for studying arbitrary rearrangements of the Haar system is
provided by the result in [J3]. There Peter Jones has given a geometric description
of homeomorphisms of the real line which preserve $\BMO$.

{\bf Acknowledgement:} I would like to thank Peter  Jones for extremely
helpful suggestions and discussions during the preparation of this paper.

\section{The geometric conditions on $\tau$}

We are now going to define Property $\cP$ and the weak-Property  $\cP$.
These definitions are difficult to unwind. We will therefore look at a
simplified version of Property $\cP$ at the end of this section.

Let $\cL$ be a collection of dyadic intervals. Following standard notation we
let $\cL\cap J=\{I\in\cL:I\sbe J\}$. $\cL^*$ denotes the set covered by
$\cL$, i.e., $\cL^*=\bigcup_{I\in\cL} I$. $\max\cL$ denotes those $K\in\cL$ that
are maximal with respect to inclusion. $\max\cL$ is a collection of pairwise
disjoint dyadic intervals which covers the same set as $\cL$.
We say that a rearrangement $\tau$ satisfies {\it Property} $\cP$ if there exists
$M>0$ so that for each $J\in\cD$ the collection $\tau(\cD)\cap J$ can be
decomposed as a disjoint union
$$\bigcup\tau(\cL_i)\cup\bigcup \cE_i$$
so that
$$\bigcup\cE_i  \mbox{ satisfies the $M$ Carleson condition,}\leqno(2.1)$$
$$\frac{|\tau(I)|}{|I|}\le M\frac {|\tau(\cL_i)^*|+|\cE_i^*|}{|\cL^*_i|},\leqno(2.2)$$
whenever $I\in \cL_i$, $i\in\tN$,
$$\sum^\infty_{i=1}|\tau(\cL_i)^*|\le M|J|.\leqno(2.3)$$

The essential, and most difficult, implication of Theorem 1 is that \newline
$||T||\, ||T^{-1}||<\infty$ causes $\tau$ and $\tau^{-1}$ to satisfy Property
$\cP$.
Later, we will describe examples of rearrangements which show that we have to
weaken Property $\cP$ to find a characterization for the case $||T||<\infty$.
We say that a rearrangement $\tau$ satisfies the {\it weak-Property} $\cP$ if
there exists $M>0$ so that for every $\cB\sb \cD$ and $J\in\cD$, the collection
$\tau(\cB)\cap J$
can be decomposed as a disjoint union
$$\bigcup\tau(\cL_i)\cup\bigcup\cE_i,$$
so that
$$\bigcup\cE_i\mbox{ satisfies the $M$-Carleson condition,}\leqno(2.4)$$
$$\frac{|\tau(I)|}{|I|}\le M\frac{|\tau(\cL_i)^*|+|\cE^*_i|}{|\cL^*_i|},\leqno(2.5)$$
whenever $I\in\cL_i$, $i\in\tN$,
$$\sum^\infty_{i=1}|\tau(\cL_i)^*|\le
M|J|\sup_i{\lbrack\!\lbrack}\tau^{-1}(\max\tau(\cL_i))
{\rbrack\!\rbrack}.\leqno (2.6)$$

The weak-Property $\cP$ allows for a nice decomposition of $\tau(\cB)\cap J$
{\it only} when $\sup_i{\lbrack\!\lbrack}\tau^{-1}(\max\tau(\cL_i))
{\rbrack\!\rbrack}$ is bounded. Hence the above condition on
$\tau$ is weaker than Property $\cP$.

It might be useful to make a few straightforward observations:
For $x=\sum x_Ih_I$ we have
$$\sup|x_I|\le||x||\leqno(2.7)$$
and for $\cL\sbe\cD$,
$$\sum_{I\in\cL}x_I^2|I|\le||x||^2|\cL^*|.\leqno(2.8)$$
Let now $x=\sum_{I\in\cL}h_I$, then $||x||={\lbrack\!\lbrack}\cL
{\rbrack\!\rbrack}^{1/2}$ and by (2.8),
$$\sum_{I\in\cL}|I|\le{\lbrack\!\lbrack}\cL{\rbrack\!\rbrack}|\cL^*|.\leqno(2.9)$$

We use the rest of this section to discuss why $\cP$ and weak $\cP$
appear {\it naturally} when boundednes of permutation operators is studied in
$\BMO$.
Suppose first, that for  a collection $\cC$ of dyadic intervals we
have
${\lbrack\!\lbrack}\tau(\cC){\rbrack\!\rbrack}\le M$. Then, in
$\BMO$, $\{h_{\tau(I)}:I\in\cC\}$ is $M^{1/2}$ equivalent to the unit vectors in
$\ell^\infty$. Indeed for $x_I\in\wR$ we have,
\begin{eqnarray*}
\left\Vert\sum_{I\in\cC}x_Ih_{\tau(I)}\right\Vert^2&=&\sup_J\frac
1{|J|}\sum_{\tau(I)\sbe J, I\in\cC}x^2_I|\tau(I)|\\
&&\\
&\le&\sup_{I\in\cC}x_I^2\sup_J\frac 1{|J|}\sum_{\tau(I)\sbe J, I\in\cC}
|\tau(I)|\\
&&\\
&\le&\sup_{I\in\cC}x_I^2M.\end{eqnarray*}
\noindent Hence, by (2.7),
$$\sup_{I\in\cC}|x_I|\le\left\Vert\sum_{I\in\cC}x_Ih_{\tau(I)}\right\Vert\le
M^{1/2}\sup_{I\in\cC}|x_I|.$$
\noindent
This shows that $\{h_{\tau(I)}:I\in\cC\}$ is equivalent to the unit vector basis
in $l^\infty$.
\noindent
By (2.7) again,
$$\sup_{I\in\cC}|x_I|\le\left\Vert\sum_{I\in\cC}x_Ih_I\right\Vert.$$
\noindent
Therefore in the case where ${\lbrack\!\lbrack}\tau(\cC){\rbrack\!\rbrack}\le
M$, the rearrangement $\tau$ does not need
to satisfy any further conditions for $T$ to be bounded on
$\span\{h_I:I\in\cC\}$.

If however ${\lbrack\!\lbrack}\tau(\cL)
{\rbrack\!\rbrack}=\infty $ for some collection $\cL$
then we should expect that a strong homogeneity condition -- to be satisfied by
$\tau$ on $\cL$ -- is necessary for $T$ to be bounded on $\span\{h_I:I\in\cL\}$.
The following condition, which is a simplified model of Property $\cP$,
specifies which homogeneity condition we have in mind.

A rearrangement $\tau$ is said to satisfy condition $\cS$, if for every
$J\in\cD$, $\tau(\cD)\cap J$ splits into $\tau(\cL)\cup \cE$ so that
$$\cE\mbox{ satisfies the $M$-Carleson condition,}\leqno(2.10)$$
and
$$\frac{|\tau(I)|} {|I|}\le
M\frac{|\tau(\cL)^*|}{|\cL^*|}, \mbox{ for } I\in\cL.\leqno(2.11)$$
Condition $\cS$ is a useful model for Property $\cP$. However, $\cS$ cannot
replace $\cP$ in Theorem 1 since $\cS$ is not a necessary condition for $T$ to
be an isomorphism on $\BMO$. (See Section 5.)
However condition $\cS$ is a sufficient condition for $T$ to be bounded. We
bring the short proof of this statement now because it illustrates
our approach to the problems addressed in this paper.

Let $x=\sum x_Ih_I$. Choose $J\in\cD$ such that
$$\frac 12||Tx||^2\le\frac 1{|J|}\sum_{\tau(I)\sbe J}x^2_I|\tau(I)|.\leqno(2.12)$$
By condition $\cS$ we may split $\{\tau(I)\sbe J\}$ as $\tau(\cL)\cup\cE$ so that
(2.10) and (2.11) hold. Write
\begin{eqnarray*}
S_1 &=&\sum_{I\in\cL}x^2_I|\tau(I)|\\
S_2&=&\sum_{\tau(I)\in\cE}x^2_I|\tau(I)|.\end{eqnarray*}
By (2.11),
$|\tau(I)|\le M|\tau(\cL)^*|\, |I|/|\cL^*|$
whenever $I\in\cL$ and so since $\tau(\cL)^*\sbe J$ and (2.8) we may estimate
$$
S_1\le M\frac{|\tau(\cL)^*|}{|\cL^*|}\sum_{I\in\cL}x^2_I|I|\leqno(2.13)$$
\begin{eqnarray*}
&\le&
M\frac{|J|}{|\cL^*|}\sum_{I\in\cL}x^2_I|I|\\
&&\\
&\le& M|J|\, ||x||^2.\end{eqnarray*}
\noindent
It follows from (2.7),  (2.9) and (2.10) that
$$S_2\le\sup x^2_I\sum_{I\in\cE}|I|
\le||x||^2\, |J|{\lbrack\!\lbrack}\cE{\rbrack\!\rbrack}
\le||x||^2\, |J|M.\leqno(2.14)$$
By (2.12), (2.13) and (2.14), $T$ is bounded on $\BMO$.

\section{Rearrangements in\-du\-cing Iso\-mor\-phi\-sms of $\BMO$}

In this section we prove the Main Lemma of this paper. It allowes us to find a
block of intervals on which $\tau$ acts homogeneously; this block will be chosen
as large as possible. We then show that the remaining intervals are comparatively
few when $||T||<\infty$. Successively applying the Main Lemma we deduce
Property $\cP$ when $||T||\,||T^{-1}||<\infty$.

For $J\in\cD$ we let
$Q(J)=\{I\in\cD,I\sbe J\}$;  for $\cL\sbe\cD$ we let
 $Q(\cL)=\{J\in\cD:\exists K\in\cL, J\sbe K\}$.

\begin{lemma} Suppose that ${\lbrack\!\lbrack}\tau(\cE){\rbrack\!\rbrack}\le M$ when $\cE$ is a collection of pairwise
disjoint dyadic intervals. Then, for $A\ge 1$ and $I_0\in\cD$ there exists
 $\cC\sb\cD$, consisting of pairwise disjoint dyadic intervals so that:

$$ \sum _{L\in\cC}|L|\le\frac{M|I_0|}A.\leqno(i)$$
(ii) If  $\cL=Q(I_0)\sm Q(\cC)$
then for  $I\in\cL$
$$\frac{|\tau(I)|}{|I|}\le A\frac{|\tau(\cL)^*|+|\tau(\cC)^*|}{|I_0|}.$$
\end{lemma}

\proof We will first construct an auxiliary collection $\cK$ of coloured dyadic
intervals. In each step of the construction we modify $\cK$ by changing the
colour of a particular $I\in\cK$, or by adding a new interval to $\cK$. We write
$\cK=\cK\cup\{L\}$ when we decide to put $L$ into the already existing
collection $\cK$.

We shall now define two construction rules and a stopping rule. Following these
rules we will then construct $\cK$ and $\cC$. Let $A\ge 1$ and pick $I_0\in\cD$.

{\it Rule 1:}  Suppose there exists a green interval $I\in\cK$ such that at
least one of the two dyadic subintervals of $I$ with length $|I|/2$ is not
contained in $\cK$. Call it $I_1$.

If
$$\frac{|\tau(I_1)|}{|I_1|}\le A\frac{|\tau(\cK\cup\{I_1\})^*|}{|I_0|},$$
then $I_1$ is a green interval and $\cK=\cK\cup\{I_1\}$.

If
$$\frac{|\tau(I_1)|}{|I_1|}> A\frac{|\tau(\cK\cup\{I_1\})^*|}{|I_0|},$$
then $I_1$ is a red interval and $\cK=\cK\cup\{I_1\}$.

{\it Rule 2:} Suppose $I\in\cK$ is a red interval.

If
$$\frac{|\tau(I)|}{|I|}\le A\frac{|\tau(\cK)^*|}{|I_0|}.\leqno(3.1)$$
Then we change the colour of $I$ from red to green. If (3.1) does not hold, then
we don't change the colour of $I$.

{\it Rule 3:} Suppose that for each red interval $I\in\cK$ we have
$$\frac{|\tau(I)|}{|I|}> A\frac{|\tau(\cK)^*|}{|I_0|},\leqno(3.2)$$
then we define $\cC$ to be the collection of red intervals in $\cK$.

Next we decree how these rules are to be applied.

The proof starts as follows. $I_0$ is a green interval and $\cK:=\{I_0\}$. Then
we apply Rule 1 until for every green interval $I$, both dyadic subintervals of
measure $|I|/2$, are contained in $\cK$.
We then apply Rule 2 to the red intervals of $\cK$.

If a new green interval is created by applying Rule 2, then we apply Rule 1, as
above, and thereafter Rule 2 again, etc.

If we don't create new green intervals by applying Rule 2 to the
red intervals of $\cK$, 
then the stopping criterion (3.2) of Rule 3 holds and we define $\cC$ to be the
red intervals in $\cK$.

Let $\cC\sb\cK\sbe\cD\cap I_0$ be the collections constructed by applying our
three rules as described above. During the construction, each interval in $\cK$
has been coloured. Note that only subintervals of green intervals can be placed
into $\cK$ by Rule 1. Hence, if $I\in\cK$ is red and $K$ is strictly contained
in $I$, then
$K\not\in\cK$. Therefore any two red intervals $L,K\in\cK$ are necessarily
disjoint. Hence, $\cC$ is a collection of pairwise disjoint intervals. By
hypothesis on $\tau$ the collection $\tau(\cC)$ satisfies the $M$ Carleson
condition. In particular, since $\tau(\cK)\spe\tau(\cC)$,
$$\frac 1{|\tau(\cK)^*|}\sum_{I\in\cC}|\tau(I)|
\le\frac 1{|\tau(\cC)^*|}\sum_{I\in\cC}|\tau(I)|\le M.\leqno(3.3)$$
On the other hand, summing the inequalities (3.2) gives
$$\frac 1{|I_0|}\sum_{I\in\cC}|I|\le\frac
1{A|\tau(\cK)^*|}\sum_{I\in\cC}|\tau(I)|.\leqno (3.4)$$

Combining (3.3) and (3.4) gives
$$\sum_{I\in\cC}|I|\le\frac{M|I_0|}A.$$
This is Lemma 1, (i).

We now turn to (ii). Recall that $K\not\in\cK$ when $K$ is strictly contained
in $I\in\cC$. Recall also that $I_0\in\cK$ and if $K\in\cK$, $I_0\spe
L\spe K$ then $L\in\cK$. If we let
$\cL=Q(I_0)\sm Q(\cC)$ 
then $\cL=\cK\sm\cC$, i.e., $\cL$ consists of the green intervals of $\cK$.
Recall that $I\in\cK$ is green when
$$\frac{|\tau(I)|}{|I|}\le A\frac{|\tau(\cK)^*|}{|I_0|}.$$
Clearly, $|\tau(\cK)^*|\le|\tau(\cL)^*|+|\tau(\cC)^*|$. Therefore,
$$\frac{|\tau(I)|}{|I|}\le A\frac{|\tau(\cL)^*|+|\tau(\cC)^*|}{|I_0|}$$
whenever $I\in\cL$. This proves Lemma 1 (ii).
\endproof

\rem The same proof shows that for any $A\ge 1$, $I_0\in\cD$ and
$\cB\sb\cD$ there exists $\cC\sb\cB$, consisting of pairwise disjoint dyadic
intervals, so that
$$\sum_{L\in\cC}|L|\le\frac{M|J|}A.\leqno{\mbox{\rm (i)}}$$
(ii) If $\cL=\cB\cap I_0\sm Q(\cC)$, 
then
$$\frac{|\tau(I)|}{|I|}\le A\frac{|\tau(\cL)^*|+|\tau(\cC)^*|}{|I_0|}.$$
\endrem

We will obtain Property $\cP$ by combining the local information obtained in
Lemma 1. The next lemma is a useful tool for achieving this.

\begin{lemma} Let $\cG_k$, $k\in\tN$ be collections of pairwise disjoint
intervals satisfying the following conditions.
$$\mbox{If } I\in\cG_k, K\in\cG_m\mbox{ and } I\cap K\ne\es,\mbox{ then }
I \sp K \mbox{ iff } k<m.\leqno(3.5)$$
$$\mbox{For } I\in\cG_k\mbox{ and } l\ge 1, \sum_{K\in\cG_{k+l}\cap I}|K|\le
2^{-l}|I|.\leqno(3.6)$$
Let $\cV_k$ be a collection of dyadic intervals in $Q(\cG_k)\sm Q(\cG_{k+1})$.
Then
$${\biggl\lbrack\!\biggl\lbrack} \bigcup\cV_k{
\biggr\rbrack\!\biggr\rbrack}\le 2\sup{\lbrack\!\lbrack}\cV_k{\rbrack\!\rbrack}.$$
\end{lemma}

\proof Let $\cV=\bigcup\cV_k$. Let $I\in\cD$. There is $k_0\in\tN$ such that
$I\in Q(\cG_{k_0})\sm Q(\cG_{k_0+1})$.
By (3.5), we may rewrite $\cV\cap I$  as  follows:
$$\cV\cap I=\bigcup^\infty_{k=k_0}\cV_k\cap I.\leqno(3.7)$$
Since $\cV_k\cap I^*$ is contained in  $\cG_k\cap I^*$, by (2.9), we obtain
$$\sum_{K\in\cV_k\cap I}|K|\le{\lbrack\!\lbrack}\cV_k{\rbrack\!\rbrack}|\cG_k\cap I^*|.\leqno(3.8)$$
By (3.6) -- (3.8), $\cV$ satisfies the
$2\sup{\lbrack\!\lbrack}\cV_k{\rbrack\!\rbrack}$ Carleson condition. For
\begin{eqnarray*}
\sum_{K\in\cV\cap I}|K|&=&\sum^\infty_{k=k_0}\sum_{K\in\cV_k\cap I}|K|\\
&=&\sum^\infty_{k=k_0}{\lbrack\!\lbrack}\cV_k{\rbrack\!\rbrack}|\cG_k\cap I^*|\\
&\le&\sup{\lbrack\!\lbrack}\cV_k{\rbrack\!\rbrack}2|I|.\end{eqnarray*}
\endproof

\begin{prop} Suppose there exists $M\ge 1$ such that ${\lbrack\!\lbrack}
\tau(\cE){\rbrack\!\rbrack}\le M{\lbrack\!\lbrack}\cE{\rbrack\!\rbrack}$ for
$\cE\sbe \cD$. Then, for $J\in\cD$, the collection $\tau(\cD)\cap J$ can be decomposed as
$$\bigcup^\infty_{i=1}\tau(\cL_i)\cup\bigcup^\infty_{i=1}\cE_i,$$
so that the following conditions hold.

\noindent (i) $\bigcup\cE_i$ satisfies the $2M$ Carleson condition.

\noindent (ii) When $I\in\cL_i$, then
$$\frac{|\tau(I)|}{|I|}\le 2M\frac{|\tau(\cL_i)^*|+|\cE_i^*|}{|\cL_i^*|}.$$


$$\sum^\infty_{i=1}|\tau(\cL_i)^*|\le|J|2M\sup_i
{\lbrack\!\lbrack}\tau^{-1}(\max\tau
(\cL_i)){\rbrack\!\rbrack}.\leqno(iii)$$
\end{prop}


\proof Let $\cG_0$ be the maximal intervals of $\tau^{-1}\{Q(J)\}$. Given
$I\in\cG_0$, let $\cC_I$ resp. $\cL_I$ satisfy the conclusion (i) resp. (ii) of
Lemma 1.
We let $\cB_I$ be the collection of dyadic intervals $K$ for which $\tilde K$
belongs to $\cC_I$. Here $\tilde K$ denotes the dyadic interval satisfying
$\tilde K\sp K$ and $|\tilde K|=2|K|$. Note that $\cB_I$ is a collection of
pairwise disjoint dyadic intervals.
Hence $\cG_1=\bigcup_I\cB_I$, where $I$ ranges over $\cG_0$, is a
collection of pairwise disjoint dyadic intervals.

We choose now $K\in\cG_1$, and we let $\cC_K$, resp. $\cL_K$, satisfy the
conclusions
(i), resp. (ii), of Lemma 1.
The collection $\cB_K$ is obtainded from $\cC_K$ by the procedure as described
for the first step.
Then $\cG_2=\bigcup_K\cB_K$, where $K$ ranges
over $\cG_1$, is a collection of pairwise disjoint dyadic intervals. Continuing
inductively we obtain a decomposition of $\tau^{-1}\{Q(J)\}$ as
$$\bigcup \cL_I\cup\bigcup\cC_I$$
where $I$ ranges over $\cG=\bigcup^\infty_{k=0}\cG_k$. Note that
$\bigcup_{I\in\cG}\cB_I$ equals $\bigcup^\infty_{k=1}\cG_k.$

We shall now verify that $\cG$ satisfies the 2 Carleson condition.
Recall that
$\cG_k$ is a collection of disjoint intervals and $\cG_{k+1}=\bigcup\cB_I$ where
$I$ ranges over $\cG_k$.
Hence, for $I\in\cG_k$, $K\in\cG_m$ with $I\cap K\ne \es$, we have
$$I\sp K\mbox{ iff } k<m.\leqno(3.9)$$
Note that by Lemma 1 (i), for $I\in\cG_k$,
$$\sum_{K\in\cG_{k+1}\cap I}|K|\le\frac{|I|}2.$$
Hence, by induction, for $I\in\cG_k$,
$$\sum_{K\in\cG_{k+l}\cap I}|K|\le\frac{|I|}{2^l}.\leqno(3.10)$$
\noindent
By (3.9) and (3.10), Lemma 2 gives
${\lbrack\!\lbrack}\cG{\rbrack\!\rbrack}\le 2$. We enumerate $\{\cC_I,I\in\cG\}$
as $\{\cC_i:i\in\tN\}$. Let $\cE_i=\tau(\cC_i)$. By hypothesis, $\tau$ preserves
the Carleson condition. Hence,
\begin{eqnarray*}
{\lbrack\!\lbrack}\bigcup\cE_i{\rbrack\!\rbrack}&=&{\lbrack\!\lbrack}\tau(\bigcup\cC_i){\rbrack\!\rbrack}\\
&\le&M{\lbrack\!\lbrack}\bigcup\cC_i{\rbrack\!\rbrack}\\
&\le&M{\lbrack\!\lbrack}\cG{\rbrack\!\rbrack}\\
&\le& M2.\end{eqnarray*}
This proves (i).

We enumerate $\{\cL_I:I\in\cG\}$ as $\{\cL_i:i\in \tN\}$ in
the same manner. Thus
Lemma 1 (ii) gives (ii). We now show (iii).
Recall that $\cL_I=Q(I)\sm Q(\cC_I)$. Hence
$$\bigcup_{I\in \cG_k}\cL_I=Q(\cG_k)\sm Q(\cG_{k+1}).\leqno(3.11)$$
Recall also that $\cL^*_I\cap\cL^*_K=\es$, when $I,K\in\cG_k$ and clearly
$\tau^{-1}(\max\tau(\cL_I))\sbe\cL_I$.
Now let $\cV_k=\bigcup 
\tau^{-1}(\max\tau(\cL_I))$,  where the union is taken over $I\in\cG_k$.
Then, by (3.9) -- (3.11), Lemma 2 gives
${\lbrack\!\lbrack}\bigcup\cV_k{\rbrack\!\rbrack}\le 2\max{\lbrack\!\lbrack}
\cV_k{\rbrack\!\rbrack}$, and so
$${\biggl\lbrack\!\biggl\lbrack}\bigcup_{I\in\cG}\tau^{-1}(\max\tau(\cL_I))
{\biggr\rbrack\!\biggr\rbrack}\le
2\sup_{I\in\cG}\lbrack\!\lbrack\tau^{-1}(\max\tau(\cL_I))\rbrack\!\rbrack.$$
Relabelling $\{\cL_I\}$ as $\{\cL_i\}$ gives (iii) for
$$\frac 1{|J|}\sum_i|\tau(\cL_i)^*|\le{\lbrack\!\lbrack} \bigcup\max\tau
(\cL_i){\rbrack\!\rbrack}\leqno(3.12)$$
\begin{eqnarray*}
&\le&M{\lbrack\!\lbrack} \bigcup\tau^{-1}(\max\tau(\cL_i)){\rbrack\!\rbrack}\\
&\le&2M\sup_i{\lbrack\!\lbrack}\tau^{-1}(\max\tau(\cL_i)){\rbrack\!\rbrack}.
\end{eqnarray*}
\endproof


\rems 1) Suppose, moreover, that $\tau^{-1}$ preserves the Carleson condition,
i.e., there exists $M\ge 1$ such that
$\lbrack\!\lbrack\tau^{-1}\cE\rbrack\!\rbrack\le
M\lbrack\!\lbrack\cE\rbrack\!\rbrack$ for $\cE\sbe\cD$. Since $\max\tau(\cL_I)$
is a collection of pairwise disjoint intervals we have
$\lbrack\!\lbrack\tau^{-1}(\max\tau(\cL_I))\rbrack\!\rbrack\le M$. Hence (3.12)
gives the estimate
$$\frac 1{|J|}\sum^\infty_{i=1}|\tau(\cL_i)^*|\le 2M^2.$$

2) The proof of Proposition 1 shows that if $\tau$ preserves the Carleson
condition, then for any $\cB\sbe\cD$, the collection $\tau(\cB)\cap J$ can be
decomposed as $\bigcup\tau(\cL_i)\cup\bigcup\cE_i$ with $\cL_i\sb\cB$ and so
that the conclusions (i) -- (iii) of Proposition 1 hold.

3) Fix $J\in\cD$ and suppose that $\tau$ preserves the Carleson condition. Let
$\cL_i$, $\cC_i$ and $\cE_i=\tau(\cC_i)$ be the collections of intervals
constructed in the proof of Proposition 1. Recall that $\cC^*_i\sbe\cL_i^*$.
Now let $\cK_i=\cL_i\cup\cC_i$, $\cE=\bigcup\cE_i$. Then by Proposition 1, we have
$$\tau(\cD)\cap J=\bigcup\tau(\cK_i),\leqno\mbox{\rm (i)}$$
$$\cE\mbox{ satisfies the $M$ Carleson condition,}\leqno\mbox{\rm (ii)}$$
$$\frac {|\tau(I)|}{|I|}\le
M\frac{|\tau(\cK_i)^*|}{|\cK_i^*|},\leqno\mbox{\rm (iii)}$$
whenever $I\in\cK_i$ and $\tau(I)\not\in\cE$,
$$\sum|\tau(\cK_i)^*|\le|J|M\sup_i\lbrack\!\lbrack\tau^{-1}
\max\tau(\cK_i)\rbrack\!\rbrack.\leqno\mbox{\rm (iv)}$$
\endrems

We now turn to the proof of Theorem 1, which has the following pattern:
\ms
\begin{tabular}{lll}
&&$\tau$ and $\tau^{-1}$ preserve the Carleson condition;\\
\phantom{xxxxxxxx}&$\Ra$&$\tau$ satisfies Property $\cP$;\\
&$\Ra$&$T$ is bounded;\\
&$\Ra$&$\tau$ preserves the Carleson condition.\end{tabular}
\ms \parindent0pt
Clearly, the hypothesis of the first line is symmetric in $\tau$ and $\tau^{-1}$.
Hence, the above statemets hold with $\tau$ replaced by $\tau^{-1}$,
respectively $T$ replaced by $T^{-1}$. The first implication is the most
difficult part of Theorem 1. It follows from the assertions of Proposition 1.
Below we prove the implication from line 2 to line 3.

\parindent0.5cm
\prooof 

(i) $\Ra$ (ii): When  $\cE\sbe\cD$ and $x=\sum_{I\in\cE}h_I$, then
$Tx=\sum_{J\in\tau(\cE)}h_J$. By (2.9), ${\lbrack\!\lbrack}\tau(\cE)
{\rbrack\!\rbrack}=||Tx||^2$
and ${\lbrack\!\lbrack}\cE{\rbrack\!\rbrack}=||x||^2$.
Thus (ii) holds with $M=||T||^2$.

(ii) $\Ra$ (iii): This part of the proof follows from Proposition 1 and Remark 1.

(iii) $\Ra$ (i): We show that $||T||<\infty$ when $\tau$ satisfies
Property $\cP$. Let $x=\sum x_Ih_I$. Then $Tx=\sum x_Ih_{\tau(I)}$. For
$J\in\cD$ we will prove that
$$\sum_{\tau(I)\sbe J}x^2_I|\tau(I)|\le 3M^2|J|\,||x||^2.\leqno(3.13)$$
Assuming $\tau$ satisfies Property $\cP$, we decompose $\tau(\cD)\cap J$ as
$$\bigcup\tau(\cL_i)\cup\bigcup\cE_i$$
so that (2.1), (2.2) and (2.3) hold.
We let $\cE=\bigcup\cE_i$ and write
\begin{eqnarray*}
S_1&=&\sum_i\sum_{I\in\cL_i} x^2_I|\tau(I)|,\\
S_2&=&\sum_{\tau(I)\in\cE} x^2_I|\tau(I)|.\end{eqnarray*}
First, we prove the estimate $S_1\le 2M^2||x||^2\,|J|$. Note that, by (2.2),
$$|\tau(I)|\le M\frac{|\tau(\cL_i)^*|+|\cE^*_i|}{|\cL_i^*|}
|I|$$
when $I\in\cL_i$.
By (2.8),
\begin{eqnarray*}
\sum_{I\in\cL_i}x^2_I|\tau(I)|&\le&
M\frac{|\tau(\cL_i)^*|+|\cE^*_i|}{|\cL^*_i|}\sum_{I\in\cL_i}x^2_I|I|\\
&\le&
M\frac{|\tau(\cL_i)^*|+|\cE^*_i|}{|\cL_i^*|}|\cL^*_i|\,||x||^2.\end{eqnarray*}
Thus, by (2.1) and (2.3),
$$S_1\le M||x||^2\sum_i|\tau(\cL_i^*)|+|\cE_i^*|\leqno(3.14)$$
\begin{eqnarray*}\phantom{S_1}&\le& M||x||^2(M+{\lbrack\!\lbrack}
\cE{\rbrack\!\rbrack})|J|\\
&&\\
&\le&2M^2||x||^2|J|.\end{eqnarray*}
Using (2.3) and (2.7) gives
$$S_2\le\sup x^2_I\sum_{\tau(I)\in\cE}|\tau(I)|\leqno(3.15)$$
$$\le M||x||^2|J|.$$
Combining (3.14) and (3.15) gives (3.13).
\endprooof

\rem Theorem 1 admits a partial extension to the case of $L^p$ spaces. There,
one considers rearrangements of the $L^p$ normalized Haar system. If $\tau$ and
$\tau^{-1}$ satisfy the Property $\cP$ then,
$$J_1:h_I\mapsto h_{\tau(I)}\leqno(3.16)$$
extends to an isomorphism on $\BMO$, and by $H^1-\BMO$ duality, the operator
$$J_0:\frac{h_I}{|I|}\mapsto\frac{h_{\tau(I)}}{|\tau(I)|}\leqno(3.17)$$
extends to an isomorphism on dyadic $H^1$.

Rearrangements of the $L^p$ normalized Haarsystem are given by the operator
$$J_{1-1/p}:\frac{h_I}{|I|^{1/p}}\mapsto\frac{h_{\tau(I)}}{|\tau(I)|^{1/p}}$$
where $1<p<\infty$. The operator $J_{1-1/p}$ coincides with the operator
$$h_I\mapsto h_{\tau(I)}\left(\frac{|\tau(I)|}{|I|}\right)^{-1/p}.$$
Hence the family of operators $\{J_t:0\le t\le 1\}$ is embedded in the {\it
analytic} family of operators
$$J_z:h_I\mapsto h_{\tau(I)}\left(\frac{|\tau(I)|}{|I|}\right)^{z-1}$$
where $z\in S=\{x+iy:0\le x\le 1, y\in\wR\}$.

>From (3.16) and (3.17) it follows that $||J_{iy}||_{H^1}\le M$ and
$||J_{1+iy}||_{\BMO}\le M$ where $M$ is independent of $y\in\wR$. Clearly
$$z\mapsto e^{-|y|}\log\left|\int T_zfg\right|$$
is uniformly bounded on $S$, whenever $f$ and $g$ are finite linear combinations
of Haar functions. Hence by complex interpolation (see [F-St]),
$||J_{1-1/p}||_{L^p}\le M$, as claimed.
We thus have shown that if $\tau$ and $\tau^{-1}$ satisfy Property $\cP$, then
permuting the $L^p$-normalized Haar system by $\tau$ leads to an isomorphism on
$L^p$.

\section{Rearrangements in\-du\-cing boun\-ded Oper\-ators on $\BMO$}

In this section we characterize rearrangements for which $||T||<\infty$; we show
that $||T||<\infty$ holds iff $\tau$ satisfies the weak-Property $\cP$.   The
pattern of the proof for Theorem 2 differs from that of Theorem 1. It is as
follows:
\ms
\begin{tabular}{lll}
\phantom{xxxx}&&$T$ is bounded;\\
&$\Ra$&$\tau$ preserves the Carleson condition;\\
&$\Ra$&$\tau$ satisfies the weak-Property $\cP$;\\
&$\Ra$&$\tau$ preserves the Carleson condition;\\
&$\Ra$&$T$ is bounded.\end{tabular}\ms\noindent
The proof of Proposition 1 gives the implication from
line 2 to line 3. For the implication from line 4 to line 5 we rephrase an
argument of P. Jones [J1, Lemma 2.1] in terms of dyadic intervals. It enters
the proof of the following proposition.

\begin{prop}  For a rearrangement $\tau$ the following are equivalent.

(i) $T:\BMO\to\BMO$ is bounded.

(ii) There exists $M\ge 1$ such that ${\lbrack\!\lbrack}\tau(\cE){\rbrack\!\rbrack}\le M{\lbrack\!\lbrack}\cE{\rbrack\!\rbrack}$ for any collection
$\cE\sb\cD$.

(iii) There exists $M\ge 1$ such that ${\lbrack\!\lbrack}\tau(\cE){\rbrack\!\rbrack}\le M$ for any collection
$\cE\sb\cD$ and ${\lbrack\!\lbrack}\cE{\rbrack\!\rbrack}\le 4$.
\end{prop}

\proof The only implication that requires proof is (iii) $\Ra$ (i). Let $x=\sum
x_Ih_I$ have norm $\le 1$. Using (iii) we prove that
$||Tx||\le M^{1/2}$.
We may assume that there exists $K\in\tN$ and $k_I\in\tN$ such that
$$x_I^2=k_IK^{-1}\mbox{ for } I\in\cD.\leqno(4.1)$$
As $x^2_I\le||x||^2\le 1$, we have $k_I\le K$.

We will  define collections of dyadic intervals $\cE_1,\dots,\cE_K$ so that
$\cE_i$ satisfies the $3$ Carleson condition, the set covered by
$\tau(\cE_i)$ is contained in $J$ and

$$\frac 1K\sum_{\tau(I)\sbe J}k_I|\tau(I)|=\frac
1K\sum^K_{i=1}\sum_{I\in\cE_i}|\tau(I)|.\leqno(4.2)$$

Let $n\in\tN$. Let $v_n$ be the vector whose entries are dyadic intervals $I$
of length $2^{-n}$ such that $\tau(I)\sbe J$. Each such $I$ appears with
multiplicity $k_I$. 

Fix $I\in\cD$ and suppose that
$|I|=2^{-n}$. We say that $I$ is in $\cE_i$ iff $I$
occupies the $p$-th position in $v_n$ and $p=i\mbox{\,\rm mod\,}K$.

Note that since  $k_I \le K$ the entries of the vectors $v_n$ are bijectively
distributed among the collections $\cE_1,\dots, \cE_K$ by the above rule. This
gives the identity (4.2). We show next that each of the $\cE_i$ satisfies the 3 Carleson
condition. Let $n\in\tN$ and $I_0\in\cD$. We take another look at
 the definition of
$\cE_i$ and see that the cardinality $A_{n,i}$ of
the collection $\{I\sbe I_0:I\in\cE_i,|I|=2^{-n}\}$ is bounded by
$$1+\frac 1K\sum_{I\sbe I_0,|I|=2^{-n}}k_I.\leqno(4.3)$$
By (4.3), (4.1) and $||x||\le 1$ we may estimate
\begin{eqnarray*} \sum_{I\in\cE_i\cap
I_0}|I|&=&\sum^\infty_{n=-\log_2|I_0|}2^{-n}A_{n,i}\\
&\le&2|I_0|+\frac 1K\sum_{I\sbe I_0}k_I|I|\\
&\le& 3|I_0|.\end{eqnarray*}
This gives ${\lbrack\!\lbrack}\cE_i{\rbrack\!\rbrack}\le 3$.
Using  identity (4.2) we now show that $||Tx||\le M^{1/2}$. Recall that the
interval $J$
was chosen so that
$$\frac 12||Tx||^2\le\frac 1{|J|}\sum_{\tau(I)\sbe J}x^2_I |\tau(I)|.$$
\noindent
By (4.1) and (4.2),
$$\sum_{\tau(I)\sbe J}x^2_I|\tau(I)|=\frac
1K\sum^K_{i=1}\sum_{I\in\cE_i}|\tau(I)|.\leqno(4.4)$$
For $i\le K$, by (2.9),
$$\sum_{I\in \cE_i}|\tau(I)|\le{\lbrack\!\lbrack}\tau(\cE_i){\rbrack\!\rbrack}|\tau(\cE_i)^*|.$$
Hence, since ${\lbrack\!\lbrack}\tau(\cE_i){\rbrack\!\rbrack}\le M$ and
$|\tau(\cE_i)^*|\le|J|$,
(4.4) is
$\le M|J|$. This proves (i) when (iii) holds.

\proooof 

(ii) $\Lra$ (i): See Proposition 2.

(ii) $\Ra$ (iii):  See Proposition 1 and the Remarks thereafter.


(iii) $\Ra$ (ii): Let $\cB\sbe\cD$ and $J\in\cD$. Assuming that (iii) holds,
i.e., assuming that
$\tau$ satisfies the weak Property $\cP$ we will first show
that
$$\sum_{\tau(I)\in \tau(\cB)\cap J}|\tau(I)|\le 3M|J|{\lbrack\!\lbrack}\cB{\rbrack\!\rbrack}^2.\leqno(4.5)$$
By (iii) we may decompose $\tau(\cB)\cap J$ as
$\bigcup\tau(\cL_i)\cup\bigcup\cE_i$ in such a way that (2.4)-(2.6) hold. Let
$S_1=\sum_i\sum_{I\in\cL_i}|\tau(I)|$ and
$S_2=\sum_i\sum_{K\in\cE_i}|K|$.
Since $\cL_i\sbe\cB$, by (2.5), (2.9) and (2.6),
\begin{eqnarray*}
S_1&\le&\sum_i\frac{|\tau(\cL_i)^*|+|\cE^*_i|}{|\cL^*_i|}\sum_{I\in\cL_i}|I|\\
&\le&{\lbrack\!\lbrack}\cB{\rbrack\!\rbrack}\sum_i|\tau(\cL_i)^*|+|\cE_i^*|\\
&\le&{\lbrack\!\lbrack}\cB{\rbrack\!\rbrack}(\lbrack\!\lbrack\cB\rbrack\!\rbrack
M+M)|J|.\end{eqnarray*}
By (2.4) we may estimate
\begin{eqnarray*} S_2&\le&\biggl{\lbrack\!\biggl\lbrack}\bigcup\cE_i
\biggr{\rbrack\!\biggr\rbrack}|J|\\
&\le&M|J|.\end{eqnarray*}
Recall that ${\lbrack\!\lbrack}\cB{\rbrack\!\rbrack}\ge 1$
whenever $\cB\ne\es$.
Combining the estimates for $S_1,S_2$ gives $S_1+S_2\le
3M{\lbrack\!\lbrack}\cB{\rbrack\!\rbrack}^2|J|$. This demostrates
(4.5). Hence
$${\lbrack\!\lbrack}\tau(\cB){\rbrack\!\rbrack}\le 3M{\lbrack\!\lbrack}\cB{\rbrack\!\rbrack}^2.
\leqno(4.6)$$
It remains to replace ${\lbrack\!\lbrack}\cB{\rbrack\!\rbrack}^2$
by $C{\lbrack\!\lbrack}\cB{\rbrack\!\rbrack}$  in (4.6).
Let $A\in\tN$ be the largest integer  which is less than
$4{\lbrack\!\lbrack}\cB{\rbrack\!\rbrack}$. By
[J2, Lemma 5.1], $\cB$
can be decomposed as $\cB_1,\dots,\cB_{A}$ in such a way
that ${\lbrack\!\lbrack}
\cB_i{\rbrack\!\rbrack}\le 4$, $i\le A$.
Applying (4.6) to $\cB_1,\dots,\cB_A$, this decomposition gives
\begin{eqnarray*}
{\lbrack\!\lbrack}\tau(\cB){\rbrack\!\rbrack}&\le&\sum^{A}_{i=1}{\lbrack\!\lbrack}\tau(\cB_i){\rbrack\!\rbrack}\\
&\le&48M{\lbrack\!\lbrack}\cB{\rbrack\!\rbrack}.\end{eqnarray*}
We have thus proved (ii).
\endproooof

\rems 1) By Theorem 2 it is clear that $T$ is an isomorphism on $\BMO$ iff
$\tau$ and $\tau^{-1}$ satisfy Property $\cP$. However, the proof of Theorem 1
given above uses only Proposition 1 and does not depend on Proposition 2.

2) It follows from the remarks after Proposition 1 that the boundednes of $T$ on
$\BMO$ can also be characterized by the following condition on $\tau$:

There exists $M\ge 1$ so that for $J\in\cD$ and $\cB\sb\cD$, there exists a
sequence $\{\cK_i\}$ of pairwise disjoint collections of intervals, and a collection
$\cE$ of intervals so that
$$\tau(\cB)\cap J=\bigcup\tau(\cK_i);\leqno\mbox{\rm (i)}$$
$$\cE\mbox{ satisfies the $M$ Carleson condition;}\leqno\mbox{\rm (ii)}$$
$$\frac{|\tau(I)|}{|I|}\le M\frac{|\tau(\cK_i)^*|}{|\cK_i^*|},\leqno
\mbox{\rm (iii)}$$
whenever $I\in\cK_i$ and $\tau(I)\not\in\cE$;
$$\sum|\tau(\cK_i)^*|\le
|J|M\sup_i\lbrack\!\lbrack\tau^{-1}\max\tau(\cK_i)\rbrack\!\rbrack.\leqno
\mbox{\rm (iv)}$$

\section{Examples}

The properties $\cP$ respectively weak-$\cP$ are clearly very similar to -- but
more complicated than -- condition $\cS$. It might be asked whether one could
use condition $\cS$ to determine when a permutation induces a bounded operator
and when not.

We describe now a rearrangement $\tau$ such that $||T||<\infty$ and $\tau$ does not
satisfy Property $\cP$.
This explains why we had to introduce the weak-Property $\cP$ in Theorem
2. The construction of $\tau$ uses a permutation $\s$ such that $||S||$
$||S^{-1}||<\infty$ but the underlying rearrangement $\s$ does not satisfy
condition $\cS$. Hence the more complicated Property $\cP$ appears in Theorem 1.

In $[0,\frac 14)$ choose pairwise disjoint dyadic intervals $K_1,\dots,K_n\dots$ and
natural numbers $l_1,\dots,l_n\dots$ recursivly so that
\begin{eqnarray*} |K_1|&=&\frac 18,\\
l_n&=& |K_n|^{-1}2^{-1},\\
|K_{n+1}|&=&|K_n|^22^{-2}.\end{eqnarray*}
The collection
$\cK_n=\{J\sb K_n:\log_2|K_n|\le\log_2|J|+l_n\}$
consists of $l_n+1$ generations:
$$G_i(\cK_n)=\{J\sb K_n:\log_2|K_n|=\log_2|J|+i\},\quad i=0,\dots,l_n.$$
$G_i(\cK_n)$ is a collection of pairwise disjoint dyadic intervals that
covers $K_n$ and
$$\sum^{l_n}_{i=1}|G_i(\cK_n)^*|=\frac 12.\leqno(5.1)$$
We define $\rho_n$ on $\cK_n$ as
$$G_i(\cK_n)\mapsto\frac 12+i|K_n|+G_i(\cK_n),$$
i.e., when $I\in G_i(\cK_n)$,  then
$$\rho_n(I)= 1/2+i|K_n|+I.$$
Notice that $\rho_n(G_i(\cK_n))^*$ is disjoint from $\rho_n(G_j(\cK_n))^*$
when
$i\ne j$.

On $\cK=\bigcup\cK_n$ we define  $\rho$ by the relation
$\rho(I)=\rho_n(I)\mbox{ iff } I\in\cK_n.$
Using the fact
that $|\rho(I)|=|I|$, we can show that $\rho$ extends  from $\cK$ to
$\cD$
injectively in such a way that the permutation operator $R$ induced by $\rho$
satisfies $||R||\le 2$.

 Notice  that, by (5.1), $\rho(\cK_n)^*=[\frac
12,1]$. Hence for $m\in\tN$,
$$\sum^m_{n=1}|\rho(\cK_n)^*|=\frac m2.\leqno(5.2)$$
Let $S_n\sb K_n$ be the dyadic interval with the same left endpoint as $K_n$
such that
$$|S_n|\, |K_n|^{-1}=\e_n\ll 1.\leqno(5.3)$$
Let $\s_n$ be the natural bijection acting between the collections
$\{J\sbe K_n:\log_2|K_n|\le\log_2|J|+l_n\}$ and
$\{J\sbe S_n:\log_2|S_n|\le\log_2|J|+l_n\}.$
On $\cK=\bigcup\cK_n$ we define $\s$ by
$\s(I)=\s_n(I) \mbox{ iff } I\in\cK_n.$
Then $\tau=\rho\circ\s^{-1}$ can be extended to $\cD$ so as to induce
a bounded operator on
$\BMO$. However, by (5.2) and (5.3), $\tau$ does not satisfy Property $\cP$
when
$\e_n\downarrow 0$ decreases sufficiently fast. Moreover, $\s$ and $\s^{-1}$ satisfy
Property $\cP$  (hence $||S||\, ||S^{-1}||<\infty$)
but $\s$ does not satisfy condition $\cS$.

\section{Transformation of Carleson measures in {\bf D} }

In this section we apply Theorem 1 to study transformations $h$ of the unit
disc
$\wD$ which preserve the class of Carleson measures. The condition on $h$
is analogous to Property $\cP$.

Let $I$ be an interval -- not necessarily dyadic -- in $[0,1)$. Then
$S(I):=\{re^{2\pi i\th}:\th\in I, 1-|I|\le r\}$. A measure $\mu$ on $\wD$ is
called a Carleson measure if, $||\mu||_C=\sup_I\mu(S(I))/|I|<\infty$, where
the
supremum is taken over all intervals in $[0,1)$.
Given $I$, we consider $T(I)=\{re^{2\pi
i\th}:\th\in I, 1-|I|<r\le 1-|I|/2\}$. If $I$ runs through all dyadic intervals
in $[0,1)$ then $\{T(I)\}$ forms a pairwise disjoint decomposition of
$\wD$.

A sequence $\{z_i\}$ in $\wD$ is called $M$-separated if, for any dyadic
interval $I$, $T(I)$ does not contain more that $M$ elements of the sequence
$\{z_i\}$. If $E=\{z_i\}$ is a $1$-separated sequence, then there exists a
uniquely determined collection $\cE=\{I_i\}$ of dyadic intervals such that
$z_i\in T(I_i)$. We let $E^*$ be the radial projection of $\{T(I_i):i\in\tN\}$
onto the unit circle $\bT$. Notice that $E^*$ naturally corresponds to $\cE^*$,
the set covered by $\cE$ in $[0,1)$ if we identify
$\bT=\{e^{2\pi i\th}:\th\in[0,1)\}$ with $[0,1)$.
We say that a transformation $h:\wD\to\wD$ satisfies Property $\cP$ if there
exists $M\ge 1$ such that every $1$-separated sequence $\{z_i\}$ can be
decomposed into $Y_1,\dots,Y_M$ so that $h\{Y_k\}$ is $1$ separated and for
every interval $J\sb[0,1)$, the sequence $h\{Y_k\}\cap S(J)$ can be decomposed
as
$$\bigcup^\infty_{i=1}h\{L_i\}\cup\bigcup^\infty_{i=1}E_i$$
so that the following conditions hold:

\noindent (i) For $E=\bigcup E_i$, the measure
$\mu=\sum_{\o\in E}(1-|\o|)\d_{\o}$
is a Carleson measure and $||\mu||_C\le M$.
$$\frac{1-|h(z)|}{1-|z|}\le M\frac{|h\{L_i\}^*|+|E_i^*|}{|L^*_i|}\leqno
\mbox{\rm (ii)}$$
whenever $z\in L_i$ and $h(z)\not\in E$.
$$\sum^\infty_{i=1}|h\{L_i\}^*|\le M|J|.\leqno\mbox{\rm(iii)}$$

The following theorem characterises those transformations of the unit disc that
preserve Carleson measures.

\begin{theor} Let $h:\wD\to\wD$ be a bijection. Then the following conditions
are equivalent:

(i) There exists $M\ge 1$ so that for every measure $\mu$ on $\wD$
$$\frac 1M||\nu||_C\le||\mu||_C\le M||\nu||_C$$
where the measure $\nu$ is defined by $\nu(A)=\mu(h^{-1}(A))$, $A\sbe\wD$

(ii) $h$ and $h^{-1}$ satisfy Property $\cP$.

(iii) There exists $M\ge 1$ so that for every sequence $\{z_i\}$ in $\wD$
$$\frac 1M||\nu||_C\le||\mu||_C\le M||\nu||_C$$
where $\mu=\sum (1-|z_i|)\d_{z_i}$
and
$\nu=\sum (1-|h(z_i)|)\d_{h(z_i)}$.
\end{theor}

The rearrangements $\tau$ described in the previous sections can be used as a
model
for the transformations $h$ appearing in Theorem 3. In fact, Theorem 1 solves a
model problem for Theorem 3. Hence, when based on Lemma 2.1 in [J1] the proof of
Theorem 3 is very similar to the proof of Theorem 1. In order to
 avoid repetition we
leave the details to the reader.
\newpage
{\bf References}\bs
[F-St] C. Fefferman, E.M. Stein, $H^p$ spaces of several variables, Acta Math.
{\bf 129} (1972), 137-193.

[J1] P.W. Jones, Carleson measures and the Fefferman-Stein decomposition of
$\BMO(I\!\!R)$, Annals Math. {\bf 111} (1980), 197-208.

[J2]  P.W. Jones, $\BMO$ and the Banach space approximation problem, Amer. J.
Math. {\bf 107} (1985), 853-893.

[J3] P.W. Jones, Homeomorphisms of the real line which preserve $\BMO$,
Ark. f. Mat. {\bf 21} (1983), 229-231.

[S1] E. Semyonov, On the equivalence in $\cL_p$ of rearrangements of the Haar
system, Soviet Math. Dokl. {\bf 19} (1978), 1292-1294.

[S2] E. Semyonov, B. St\"ockert, Haar system rearrangement in the spaces
$\cL_p$, Analysis Math. {\bf 7} (1981), 277-298.

[S3] E. Semyonov, On the boundedness of operators rearranging the Haar system
in
the space $\cL_p$, Int. Eq. Op. Th. {\bf 6} (1983), 385-404.
\bs
\parindent0pt
Current address:\\
Yale University, 06520 New Haven, CT\\
e-mail: muller@math.yale.edu

\end{document}